\documentclass{article}
\newtheorem{lema}{Lemma}
\newtheorem{prop}{Proposition}
\newtheorem{coro}{Corolary}
\newtheorem{defi}{Definition}

\newtheorem{ejem}{Example}
\newtheorem{nota}{Note}
\usepackage{amsmath}
\usepackage{amssymb}

\usepackage[english]{babel}
\usepackage{float}
\usepackage{hyperref}
\usepackage{graphicx} 
\title{Visually-friendly manifolds with arbitrary finite fundamental group}

\author{Luca Tanganelli Castrillón}
\date{June 2024}

\begin{document}
\maketitle

\begin{abstract}
We exhibit a family of metrizable manifolds such that any finite group appears as the fundamental group of one of them. These spaces are especially interesting as they can be easily visualized, as opposed to classical examples of spaces with arbitrary fundamental group.
\end{abstract}

\section{Metrizing the quotient}
Let $G$ be a group acting on a metric space $(X,d)$. We say that $G$ acts isometrically on $X$ if for each $g\in G$ the map $x\mapsto gx$ is an isometry. We say that $G$ acts closedly on $X$ if the orbit of any element of $X$ is a closed subset of $X$. Note that it suffices for $G$ to be finite for the action to be closed.

Moreover we recall the modulo relation that $G$ induces on $X$: $x\sim y\iff\exists g\in G:x=gy\iff \mathcal{O}_x=\mathcal{O}_y$ where $\mathcal{O}_x$ is the orbit of $x$. The quotient space is denoted $X/G$, with canonical projection $p:X\to X/G$. The equivalence class of $x\in X$ is denoted $\overline x$.

\begin{prop}
    Let $G$ be a group acting isometrically and closedly on a metric space $(X,d)$. Then $X/G$ inherits a natural metric $\overline d$ given by $$\overline d(\overline x, \overline y)=d(x,\mathcal{O}_y) =d(\mathcal{O}_x,\mathcal{O}_y) =\min_{g\in G}d(x,gy)$$
    and $\overline d$ induces the quotient topology on $X/G$.
\end{prop}
\textbf{Proof.} It's easy to check that $\overline d$ is indeed a metric, and that for all $\epsilon >0$ and $x\in X$,
$$p^{-1}(B_{\overline d}(\overline x; \epsilon))=\bigcup_{g\in G}B_d(gx;\epsilon)=\bigcup_{g\in G} gB_d(x;\epsilon).$$
To see that $\overline d$ induces the quotient topology, let $A$ be an open subset of $X/G$. Let $\overline x\in A$ and let $\epsilon$ such that $B_d(x;\epsilon)\subset p^{-1}(A)$ which is open in $X$. Since $p^{-1}(A)$ is stable by $G$, it follows that $\bigcup_{g\in G} gB_d(x;\epsilon)\subset p^{-1}(A)$ and we're done.
\begin{nota}
     If $(X,d)$ is complete, then so is $(X/G, \overline d)$.
\end{nota}
\textbf{Proof.} Let $(x_n)_{n=1}^\infty$ be a sequence in $X$ such that the sequence $(\overline{x_n})$ is Cauchy. By taking a suitable subsequence, we may suppose without loss of generality that for all $n\geq 1$, $\overline{d}(\overline{x_n},\overline{x_{n+1}})\leq 2^{-n}$. Thus for each $n\geq 1$, we can choose some $g_n\in G$ such that $d(x_n,g_nx_{n+1})\leq 2^{-n}$. Now define a new sequence $(y_n)$ by $$y_1=x_1,\quad y_n=g_1g_2\ldots g_{n-1}x_n,\,n\geq 2.$$
Then $d(y_n,y_{n+1})=d(g_1\ldots g_{n-1}x_n,g_1\ldots g_nx_{n+1})=d(x_n,g_nx_{n+1})\leq 2^{-n}$. It's inmediate that $(y_n)$ is a Cauchy sequence and thus converges. By the continuity of $p$, $(\overline{y_n})$ converges. And $\overline{y_n}=\overline{x_n}$ so we're done. 

\begin{ejem}
    $\Sigma_n$ acts isometrically on $\mathbb{C}^n$ by permuting the coordinates of each vector, and $\mathbb{C}^n/\Sigma_n\cong \mathbb{C}^n$.
\end{ejem}
\textbf{Proof.} The map $\phi:\mathbb{C}^n\to\mathbb{C}^n$ given by $\phi(z_1,\ldots,z_n)=\text{coef}((z-z_1)\ldots (z-z_n))$ is continuous and compatible with the relation modulo $\Sigma_n$. (The function $\text{coef}$ extracts the $n$ ordered non-trivial coefficients of the monic polynomial that is passed as an argument.) Thus the induced map $\Tilde{\phi}:\mathbb{C}^n/\Sigma_n\to\mathbb{C}^n$ is continuous and of course bijective. In particular, it is an embedding when restricted to any compact subset of $\mathbb{C}^n/\Sigma_n$ (*). To see that $\Tilde{\phi}^{-1}$ is continuous, let $(\mathbf{z}_k)_{k=1}^\infty\subset \mathbb{C}^n$, where $\mathbf{z}_k=(z^1_k,\ldots,z^n_k)$, be a sequence such that the sequence $(\Tilde{\phi}(\overline{\mathbf{z}_k}))$ converges. We wish to see that $(\overline{\mathbf{z}_k})$ converges.

Since $(\Tilde{\phi}(\overline{\mathbf{z}_k}))$ is bounded, the polynomials $(z-z^1_k)\ldots(z-z^n_k)$ have bounded coefficients and, as is a standard result, bounded roots. Thus $(\mathbf{z}_k)$ is bounded, say contained in some closed ball $\bar B$. This implies that $(\overline{\mathbf{z}_k})$ is contained in $p(\bar B)$, which is compact, and by (*) it must converge as well.

\section{The fundamental group of the quotient}

Let $G$ be a group acting on a metric space $(X,d)$. We say $G$ acts discretely on $X$ if for all $x\in X$ the orbit of $x$ is a discrete subspace of $X$. We say $G$ acts injectively on $X$ if, for all $g,h\in G$ and $x\in X$, $gx=hx\implies g=h$ (equivalently, only the identity in $G$ can fix points of $X$).

\begin{prop}
    Let $G$ be a group acting isometrically, closedly, discretely and injectively on a metric space $(X,d)$. Then the projection $p:X\to X/G$ is a covering space.
\end{prop}
\textbf{Proof.} Let $x\in X$. Let $\epsilon=\inf_{g\in G\setminus\{1\}}d(x,gx)>0$. Let $U=B_{\overline d}(\overline x; \delta)$, where $\delta=\frac\epsilon 5$. We then observe that
$$p^{-1}(U)=\bigsqcup_{g\in G}B(gx;\delta).$$
For each $g\in G$, let $B_g=B(gx;\delta)$. We will see that $p|_{B_g}$ is an embedding.

(1) Injectivity: let $w,z\in B_g$ such that $\overline w=\overline z$. Then there is some $h\in G$ such that $w=hz$. Hence
$$d(x,g^{-1}hgx)=d(gx,hgx)\leq d(gx,hz)+d(hz,hgx)=d(gx,w)+d(z,gx)<2\delta<\epsilon$$
and so $g^{-1}hg=1\implies h=1\implies w=z$.

(2) Surjectivity: let $\overline y\in U$ and $h\in G$ such that $y\in B_h$. Then
$$d(hx,y)<\delta \implies d(gx,gh^{-1}y)<\delta\implies gh^{-1}y\in B_g.$$

(3) Homeomorphism onto its image. It will suffice to prove that for all $w,z\in B_g$, $d(z,w)=\overline d(\overline z, \overline w)$. Let $h\in G$ such that $d(z,hw)\leq d(z,w)$. We want to show $h=1$. We have
\begin{align*}
    d(x,g^{-1}hgx)=d(gx,hgx)&\leq d(gx,z)+d(z,hw)+d(hw,hgx)\\
                            &\leq d(gx,z)+d(z,w)+d(w,gx)\\
                            &\leq d(gx,z)+d(z,gx)+d(gx,w)+d(w,gx)<4\delta<\epsilon
\end{align*}
hence $g^{-1}hg=1$ and $h=1$. Thus $\overline d(\overline z, \overline w)=\min_{h\in G}d(z,hw)=d(z,w)$.

\begin{prop}
    Under the conditions stated above, if $X$ is simply connected, then $\pi_1(X/G)\cong G$.
\end{prop}
\textbf{Proof.} Let $x_0\in X$. For each $g\in G$, we let $\sigma_g:[0,1]\to X$ be any path such that $\sigma_g(0)=x_0$ and $\sigma_g(1)=gx_0$. Now define $\phi:G\to \pi_1(X/G, \overline{x_0})$ by
$$\phi(g)=[p\circ \sigma_g].$$

(1) $\phi$ is a homomorphism. Let $g,h\in G$. Then observe that $\sigma_{gh}\sim \sigma_g\cdot (g\sigma_h)$ rel $\{0,1\}$. Thus $p\circ \sigma_{gh}\sim (p\circ \sigma_g)\cdot(p\circ (g\sigma_h))$ rel $\{0,1\}$. But $p\circ (g\sigma_h)=p\circ \sigma_h$ and so
$$\phi(gh)=[p\circ\sigma_{gh}]=[(p\circ \sigma_g)\cdot(p\circ \sigma_h)]=[(p\circ \sigma_g)][(p\circ \sigma_h)]=\phi(g)\phi(h).$$

(2) $\phi$ is surjective. Let $\alpha \in \Omega(X/G,\overline{x_0})$. Since $p:X\to X/G$ is a covering space, we can lift $\alpha$ to some $\Tilde{\alpha}_{x_0}:[0,1]\to X$ with basepoint $x_0$. Now, we know $\Tilde{\alpha}_{x_0}(1)=gx_0$ for some $g\in G$, and, since $X$ is simply connected, it follows that $\Tilde{\alpha}_{x_0}\sim\sigma_g$ rel $\{0,1\}$, and so $\phi(g)=[p\circ \sigma_g]=[p\circ \Tilde{\alpha}_{x_0}]=[\alpha]$.

(3) $\phi$ is injective. Let $g,h\in G$ such that $\phi(g)=\phi(h)$. Then $p\circ\sigma_g\sim p\circ\sigma_h$ rel $\{0,1\}$. Now observe that $\sigma_g$ and $\sigma_h$ are the liftings of $p\circ\sigma_g$ and $ p\circ\sigma_h$, respectively,  based at $x_0$. By the properties of covering spaces, we know that $\sigma_g\sim\sigma_h$ rel $\{0,1\}$. Therefore $gx_0=\sigma_g(1)=\sigma_h(1)=hx_0$. By the injectiveness of the action of $G$, it follows that $g=h$.

\section{One-for-all}
Consider $E=\mathbb{R}^3$ with the usual topology. We can consider the space $X^n=\{(x_1,\ldots,x_n)\in E^n:i\neq j\implies x_i\neq x_j\}$, the space of all tuples of $n$ pairwise distinct points of $E$.

Let $G$ be any subgroup of $\Sigma_n$. Then $G$ acts on $X^n$ by permuting the components of each point: for all $\sigma \in G$, $(x_1,\ldots,x_n)\in X^n$,
$$\sigma*(x_1,\ldots,x_n)=(x_{\sigma(1)},\ldots x_{\sigma(n)}).$$
This action is isometric, closed, discrete and injective, the latter being precisely because the elements of $X$ have pairwise different coordinates. We'd now like to show that $X^n$ is simply connected.

Note that $X^n=E^n\setminus \bigcup_{1\leq i<j\leq n}L_{ij}$ where $L_{ij}=\{(x_1,\ldots,x_n)\in E^n:x_i=x_j\}$. Each $L_{ij}$ is a $(3(n-1))$-dimensional hyperplane. This motivates the following two lemmas.
\begin{lema}
    Let $\alpha\subset \mathbb{R}^n$ be an $(n-3)$-dimesional hyperplane. Let $P_1,P_2,P_3\in\mathbb{R}^n$ be non collinear, with $U_i$ a neighbourhood of $P_i$, $i=1,2,3$. Then there exist open sets $V_i\subset U_i$, $i=1,2,3$, such that for any choice of $Q_i\in V_i$, $i=1,2,3$, the triangle $\triangle (Q_1Q_2Q_3)$ does not intersect $\alpha$.
\end{lema}
\textbf{Proof.} Let $\tau=\triangle(P_1P_2P_3)$. If $\tau\cap \alpha=\emptyset$, we're done as $d(\tau,\alpha)>0$ and $d(\triangle(XYZ),\alpha)$ varies continuously in $(X,Y,Z)$ around $(P_1,P_2,P_3)$.

Suppose there is $P\in \tau\cap\alpha$. We have the inclusion $\tau \subset P+V$ where $V=\langle P_2-P_1,P_3-P_1\rangle$. Let $W$ be the directing subspace of $\alpha$. Let $t\in V^\perp\cap W^\perp$ be nonzero. We'll see that for any $\lambda>0$, $d(\tau+\lambda t, A)>0$ and thus $\triangle((P_1+\lambda t)(P_2+\lambda t)(P_3+\lambda t))\cap \alpha =0$. By choosing a small enough $\lambda$ and arguing same as in the first sentence of the proof, we'll be done

Let $Q=P+v,v\in V$. Let $S=\text{proj}_W(Q)=P+\text{proj}_W(v)$. Then naturally $\text{proj}_W(Q+\lambda t)=S$. Moreover $t\cdot (Q-S)=t\cdot(Q-P)+t\cdot (P-S)=0$
so
$$d(Q+\lambda t,\alpha)^2=d(Q+\lambda t,S)^2=\lambda^2||t||^2+d(Q,S)^2>0.$$
\begin{coro}
    Let $\alpha_1,\ldots,\alpha_k\subset \mathbb{R}^n$ be $(n-3)$-dimesional hyperplanes. Let $P_1,P_2,P_3\in\mathbb{R}^n$ be non collinear, with $U_i$ a neighbourhood of $P_i$, $i=1,2,3$. Then there exist points $Q_i\in U_i$, such that $\triangle (Q_1Q_2Q_3)\cap\alpha_j=\emptyset$ for all $j=1,\ldots,k$.
\end{coro}
\textbf{Proof.} Apply the previous lemma $k$-times, once for each hyperplane.
\begin{defi}
    A path is called a segment if its graph is a segment.
\end{defi}
\begin{lema}
    Let $A$ be a locally convex topological subspace of $\mathbb{R}^n$, and $\sigma:[0,1]\to A$ a path. Then $\sigma$ is homotopic relative to $\{0,1\}$ to a finite chain of segments.
\end{lema}
\textbf{Proof.} For each $t\in [0,1]$, let $V_t$ be an convex neighbourhood of $\sigma(t)$ and let $U_t$ be an open, connected neighbourhood of $t$ such that $\sigma(U_t)\subset V_t$. By compacity of the unit interval, we can write $[0,1]=U_{t_1}\cup\ldots\cup U_{t_s}$. Perform some homotopies in each $U_{t_i}$ and the result comes out inmediately.

\begin{figure}[H]
    \includegraphics[width=01\textwidth]{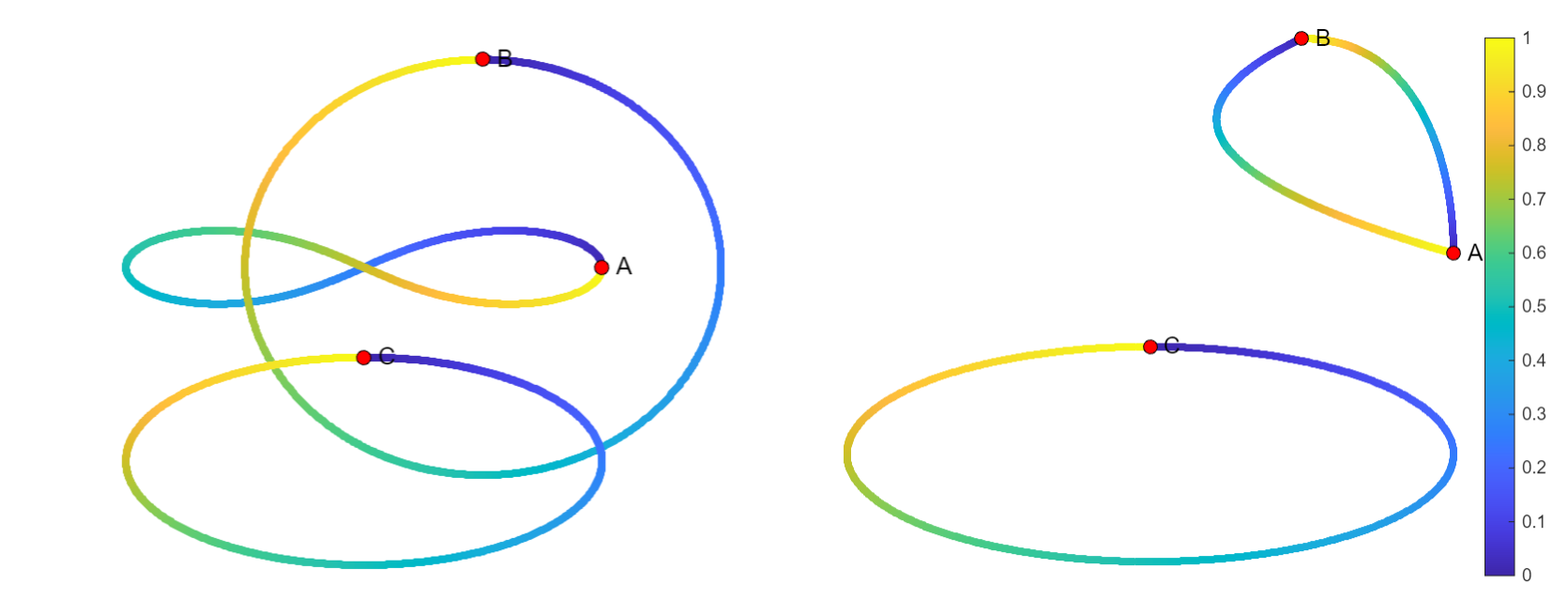}
    \caption{Two examples of loop based at $\overline{(A,B,C)}\in X^3/\Sigma_3$. The first one is null-homotopic. The second isn't because its lift to $X^3$ ``permutes $A$ and $B$''.}
    \label{fig:enter-label}
\end{figure}
\begin{prop}
    Let $\alpha_1,\ldots,\alpha_k\subset \mathbb{R}^n$ be $(n-3)$-dimesional hyperplanes. Then $\mathbb{R}^n\setminus (\alpha_1\cup\ldots\cup\alpha_k)$ is simply connected. In particular $X^n$ is simply connected.
\end{prop}
\textbf{Proof.} Let $R=\alpha_1\cup\ldots\cup\alpha_k$ and $Y=\mathbb{R}^n\setminus R$. $Y$ is an open subset of $\mathbb{R}^n$, hence locally convex. Let $\sigma:[0,1]\to Y$ be a loop, $Q_0=\sigma(0)$. By Lemma 2, $\sigma$ is homotopic to some loop
$$s_{Q_0Q_1} s_{Q_1Q_2}\ldots s_{Q_mQ_0}$$
where $Q_1,\ldots,Q_m\in Y$ and $s_{Q_iQ_j}$ is a segment from $Q_i$ to $Q_j$. WLOG we can suppose that three consecutive $Q_i$'s are never collinear (otrherwise join their segments). For $T_1,\ldots,T_s\in\mathbb{R}^n$, define the polygon $F(T_1,\ldots,T_s)=[T_1T_2]\cup\ldots\cup[T_{s-1}T_s]$. Then $d(F(\cdot),R)$ is continuous. In particular, since $d(F(Q_0,Q_1,\ldots,P),R)>0$, there exist $U_i\subset Y$ neighbourhoods of $Q_i$, for $i=0,1,2$, such that for any choice of $P_i\in U_i$, we have $d(F(P_0,P_1,P_2,Q_3,\ldots,Q_m,P_0),R)>0$. By Corollary 1 we can take $P_i\in U_i$ such that $\triangle(P_0P_1P_2)\cap R=\emptyset$. Thus by moving linearly $Q_i$ to $P_i$ it's clear that
$$\sigma\sim s_{P_0P_1}s_{P_1P_2}s_{P_2Q_3}\ldots s_{Q_mP_0}$$
and since $s_{P_0P_1}s_{P_1P_2}\sim s_{P_0P_2}$, we get
$$\sigma\sim s_{P_0P_2}s_{P_2Q_3}\ldots s_{Q_mP_0}.$$
We have decremented by one the number of segments in the chain. Continue the process until having a relation of the sort $\sigma\sim s_{AB}s_{BA}$ which implies that $\sigma$ is null-homotopic.

\begin{coro}
    $\pi_1(X^n/G)\cong G$.
\end{coro}
\begin{prop}
    $X^n/G$ is a metrizable $(3n)$-dimensional manifold.
\end{prop}
\textbf{Proof.} For the local resemblance to $\mathbb{R}^{3n}$, see proof of Proposition 2.
\begin{ejem}
    $X^2/\Sigma_2\cong \mathbb{R}^4\times \mathbb{R}P^2$.
\end{ejem}
A loop in $X^n/G$ can be visualized as $n$ ``simultaneous'' paths in $\mathbb{R}^3$ that never meet and such that the tuple of ending points is a permutation (from $G$) of the tuple of starting points.

In practice, these representations can be simply drawn in $2$-dimensional space, allowing some points to momentarily ``jump'' on other points when performing homotopies in order to avoid collisions, similar to how knots are drawn on a sheet of paper.

\end{document}